\documentclass[english]{amsart}
\usepackage[T1]{fontenc}
\usepackage[latin9]{inputenc}
\usepackage{babel}
\usepackage{textcomp}
\usepackage{amstext}
\usepackage{amsthm}
\usepackage[unicode=true,pdfusetitle,
 bookmarks=true,bookmarksnumbered=false,bookmarksopen=false,
 breaklinks=false,pdfborder={0 0 1},backref=false,colorlinks=false]
 {hyperref}

\makeatletter
\numberwithin{equation}{section}
\theoremstyle{plain}
\newtheorem{thm}{\protect\theoremname}[section]
\theoremstyle{plain}
\newtheorem{cor}{\protect\corollaryname}[section]

\makeatother

\providecommand{\corollaryname}{Corollary}
\providecommand{\theoremname}{Theorem}

\begin{document}
\title[$q$-Trigonometric identities]{On certain $q$-trigonometric identities  }
\author{Bing He}
\address{School of Mathematics and Statistics, Central South University\\
 Changsha 410083, Hunan, People's Republic of China}
\email{yuhe001@foxmail.com; yuhelingyun@foxmail.com}
\begin{abstract}
Finding theta function (or $q$-)analogues for well-known trigonometric
identities is an interesting topic. In this paper, we first introduce
the definition of $q$-analogues for $\mathrm{tan}z$ and $\mathrm{cot}z$
and then apply the theory of elliptic functions to establish a theta
function identity. From this identity we deduce two $q$-trigonometric
identities involving $\mathrm{tan}_{q}z$ and $\cot_{q}z,$ which
are theta function analogues for two well-known trigonometric identities
concerning $\mathrm{tan}z$ and $\cot z.$ Some other $q$-trigonometric
identities are also given.
\end{abstract}

\maketitle

\section{Introduction}

To carry out our study, we need the definition of Jacobi's theta functions
\cite{GR,L01,WW}:
\begin{align*}
\theta_{1}(z|\tau) & =-iq^{\frac{1}{4}}\sum_{k=-\infty}^{\infty}(-1)^{k}q^{k(k+1)}e^{(2k+1)zi},\;\theta_{3}(z|\tau)=\sum_{k=-\infty}^{\infty}q^{k^{2}}e^{2kzi},\\
\theta_{2}(z|\tau) & =q^{\frac{1}{4}}\sum_{k=-\infty}^{\infty}q^{k(k+1)}e^{(2k+1)zi},\;\theta_{4}(z|\tau)=\sum_{k=-\infty}^{\infty}(-1)^{k}q^{k^{2}}e^{2kzi},
\end{align*}
where $q=\exp(\pi i\tau)$ and $\text{Im}\:\tau>0.$ 

Throughout this paper, we use the notation $\vartheta_{j}(\tau)$
to denote $\theta_{j}(0|\tau)$ for $j=2,3,4$. 

The Jacobi infinite product expressions for the theta functions are
well-known:
\[
\begin{aligned} & \theta_{1}(z|\tau)=2q^{\frac{1}{4}}\sin z(q^{2};q^{2})_{\infty}(q^{2}e^{2zi};q^{2})_{\infty}(q^{2}e^{-2zi};q^{2})_{\infty},\\
 & \theta_{2}(z|\tau)=2q^{\frac{1}{4}}\cos z(q^{2};q^{2})_{\infty}(-q^{2}e^{2zi};q^{2})_{\infty}(-q^{2}e^{-2zi};q^{2})_{\infty},\\
 & \theta_{3}(z|\tau)=(q^{2};q^{2})_{\infty}(-qe^{2zi};q^{2})_{\infty}(-qe^{-2zi};q^{2})_{\infty},\\
 & \theta_{4}(z|\tau)=(q^{2};q^{2})_{\infty}(qe^{2zi};q^{2})_{\infty}(qe^{-2zi};q^{2})_{\infty},
\end{aligned}
\]
where $(z;q)_{\infty}$ is the $q$-shifted factorial given by
\[
(z;q)_{\infty}=\prod_{n\geq0}(1-zq^{n}).
\]
With respect to the (quasi) periods $\pi$ and $\pi\tau,$ we have
\begin{equation}
\begin{gathered}\theta_{1}(z+\pi|\tau)=\text{\textminus}\theta_{1}(z|\tau),\quad\theta_{1}(z+\pi\tau|\tau)=-q^{\text{\textminus}1}e^{\text{\textminus}2zi}\theta_{1}(z|\tau),\\
\theta_{2}(z+\pi|\tau)=\text{\textminus}\theta_{2}(z|\tau),\quad\theta_{2}(z+\pi\tau|\tau)=q^{\text{\textminus}1}e^{\text{\textminus}2zi}\theta_{2}(z|\tau),\\
\theta_{3}(z+\pi|\tau)=\theta_{3}(z|\tau),\quad\theta_{3}(z+\pi\tau|\tau)=q^{\text{\textminus}1}e^{\text{\textminus}2zi}\theta_{3}(z|\tau),\\
\theta_{4}(z+\pi|\tau)=\theta_{4}(z|\tau),\quad\theta_{4}(z+\pi\tau|\tau)=-q^{\text{\textminus}1}e^{\text{\textminus}2zi}\theta_{4}(z|\tau).
\end{gathered}
\label{eq:1-4}
\end{equation}
For the half period $\pi\tau/2,$ we also have
\begin{equation}
\begin{gathered}\theta_{1}\bigg(z+\frac{\pi\tau}{2}\bigg|\tau\bigg)=iB\theta_{4}(z|\tau),\quad\theta_{2}\bigg(z+\frac{\pi\tau}{2}\bigg|\tau\bigg)=B\theta_{3}(z|\tau),\\
\theta_{3}\bigg(z+\frac{\pi\tau}{2}\bigg|\tau\bigg)=B\theta_{2}(z|\tau),\quad\theta_{4}\bigg(z+\frac{\pi\tau}{2}\bigg|\tau\bigg)=iB\theta_{1}(z|\tau),
\end{gathered}
\label{eq:1-3}
\end{equation}
where $B=q^{-1/4}e^{-iz}.$ 

In \cite{G} Gosper introduced $q$-analogues of $\sin z$ and $\cos z:$
\begin{eqnarray*}
 &  & \sin_{q}(\pi z):=\frac{(q^{2-2z};q^{2})_{\infty}(q^{2z};q^{2})_{\infty}}{(q;q^{2})_{\infty}^{2}}q^{(z-1\text{/2})^{2}},\\
 &  & \cos_{q}(\pi z):=\frac{(q^{1-2z};q^{2})_{\infty}(q^{1+2z};q^{2})_{\infty}}{(q;q^{2})_{\infty}^{2}}q^{z^{2}},
\end{eqnarray*}
and gave two relations between $\sin_{q},\;\cos_{q}$ and the functions
$\theta_{1}$ and $\theta_{2},$ which are equivalent to 
\[
\sin_{q}z=\frac{\theta_{1}(z|\tau')}{\vartheta_{2}(\tau')},\;\cos_{q}z=\frac{\theta_{2}(z|\tau')}{\vartheta_{2}(\tau')},
\]
where $\tau'=-\frac{1}{\tau}.$ It is easily seen that $\cos_{q}z=\sin_{q}(\frac{\pi}{2}\pm z).$
In that paper, using empirical evidence based on a computer program
called MACSYMA, Gosper conjectured without proofs various identities
involving $\sin_{q}z$ and $\cos_{q}z.$ Some of these conjectures
were confirmed by different authors \cite{AAE,E1,E2,M}. 

We now define the $q$-analogues for $\mathrm{tan}z$ and $\mathrm{cot}z:$
\[
\mathrm{tan}_{q}(\pi z):=\frac{\sin_{q}(\pi z)}{\cos_{q}(\pi z)}=\dfrac{(q^{2-2z};q^{2})_{\infty}(q^{2z};q^{2})_{\infty}}{(q^{1-2z};q^{2})_{\infty}(q^{1+2z};q^{2})_{\infty}}q^{1/4-z}
\]
and 
\[
\mathrm{cot}_{q}(\pi z):=\frac{\cos_{q}(\pi z)}{\sin_{q}(\pi z)}=\dfrac{(q^{1-2z};q^{2})_{\infty}(q^{1+2z};q^{2})_{\infty}}{(q^{2-2z};q^{2})_{\infty}(q^{2z};q^{2})_{\infty}}q^{z-1/4}.
\]
Then
\begin{equation}
\mathrm{tan}_{q}z=\frac{\theta_{1}(z|\tau')}{\theta_{2}(z|\tau')}\label{eq:0-2}
\end{equation}
and 
\[
\mathrm{cot}_{q}z=\frac{\theta_{2}(z|\tau')}{\theta_{1}(z|\tau')}.
\]

Trigonometric identities is a very importan topic. Two well-known
trigonometric identity are as follows \cite{W}: if $x+y+z=\pi,$
then 
\begin{align}
 & \mathrm{tan}x+\mathrm{tan}y+\mathrm{tan}z=\mathrm{tan}x\mathrm{tan}y\mathrm{tan}z,\label{eq:0-1}\\
 & \mathrm{cot}x\mathrm{cot}y+\mathrm{cot}y\mathrm{cot}z+\mathrm{cot}z\mathrm{cot}x=1.\label{eq:00-1}
\end{align}
Finding theta function analogues for the trigonometric identities
\eqref{eq:0-1} and \eqref{eq:00-1} is also interesting. Our motivation
for the present work emanates from \cite{G}. In this paper we shall
establish the following $q$-trigonometric identities. 
\begin{thm}
\label{t1} If $x+y+z=\pi,$ then
\begin{equation}
\begin{gathered}\mathrm{ccs}_{q}(x-y)\mathrm{tan}_{q^{2}}x+\mathrm{ccs}_{q}(x-y)\mathrm{tan}_{q^{2}}y+\mathrm{ssn}_{q}(x-y)\mathrm{tan}_{q}z\\
=\mathrm{ssn}_{q}(x-y)\mathrm{tan}_{q^{2}}x\mathrm{tan}_{q^{2}}y\mathrm{tan}_{q}z,
\end{gathered}
\label{eq:1-6}
\end{equation}
and 
\begin{equation}
\mathrm{ssn}_{q}(x-y)\mathrm{cot}_{q^{2}}x\mathrm{cot}_{q^{2}}y+\mathrm{ccs}_{q}(x-y)\mathrm{cot}_{q^{2}}y\mathrm{cot}_{q}z+\mathrm{ccs}_{q}(x-y)\mathrm{cot}_{q}z\mathrm{cot}_{q^{2}}x=\mathrm{ssn}_{q}(x-y),\label{eq:1-7}
\end{equation}
where $\mathrm{ccs}_{q}z$ and $\mathrm{ssn}_{q}z$ are also $q$-trigonometric
functions given by 
\[
\mathrm{ccs}_{q}z:=\frac{\cos_{q^{2}}z}{\cos_{q}z},\;\mathrm{ssn}_{q}z:=\frac{\sin{}_{q^{2}}z}{\sin_{q}z}.
\]
\end{thm}
Setting $q\rightarrow1$ in \eqref{eq:1-6} and \eqref{eq:1-7} and
noting that 
\[
\lim_{q\rightarrow1}\big\{\mathrm{tan}_{q}z,\mathrm{ssn}_{q}z,\mathrm{ccs}_{q}z\big\}=\big\{\mathrm{tan}z,1,1\big\},
\]
we can easily obtain \eqref{eq:0-1} and \eqref{eq:00-1} respectively. 

We replace $(x,y,z)$ by $(\pi/2-x,\pi/2-y,\pi/2-z)$ in \eqref{eq:1-6}
and \eqref{eq:1-7} to get 
\begin{cor}
If $x+y+z=\pi/2,$ then
\[
\begin{gathered}\mathrm{ccs}_{q}(x-y)\mathrm{cot}_{q^{2}}x+\mathrm{ccs}_{q}(x-y)\mathrm{cot}_{q^{2}}y+\mathrm{ssn}_{q}(x-y)\mathrm{cot}_{q}z\\
=\mathrm{ssn}_{q}(x-y)\mathrm{cot}_{q^{2}}x\mathrm{cot}_{q^{2}}y\mathrm{cot}_{q}z
\end{gathered}
\]
and 
\[
\mathrm{ssn}_{q}(x-y)\mathrm{tan}_{q^{2}}x\mathrm{tan}_{q^{2}}y+\mathrm{ccs}_{q}(x-y)\mathrm{tan}_{q^{2}}y\mathrm{tan}_{q}z+\mathrm{ccs}_{q}(x-y)\mathrm{tan}_{q}z\mathrm{tan}_{q^{2}}x=1.
\]
\end{cor}
These identities are $q$-analogues of the trigonometric identities:
if $x+y+z=\pi/2,$ then
\begin{align*}
 & \mathrm{cot}x+\mathrm{cot}y+\mathrm{cot}z=\mathrm{cot}x\mathrm{cot}y\mathrm{cot}z,\\
 & \mathrm{tan}x\mathrm{tan}y+\mathrm{tan}y\mathrm{tan}z+\mathrm{tan}z\mathrm{tan}x=1.
\end{align*}

In order to prove Theorem \ref{t1} we need to establish the following
theta function identity by employing the theory of elliptic functions.
For more information dealing with formulas of theta functions by using
elliptic functions, please, see Liu \cite{L05,L5,L09,L07} and Shen
\cite{S1,S2}.
\begin{thm}
\label{t2} For all complex numbers $x$ and $y,$ we have 
\begin{equation}
\begin{gathered}\theta_{2}(x+y|2\tau)\theta_{3}(x-y|2\tau)(\theta_{1}(x|\tau)\theta_{2}(y|\tau)+\theta_{1}(y|\tau)\theta_{2}(x|\tau))\\
=\theta_{1}(x+y|2\tau)\theta_{4}(x-y|2\tau)(\theta_{2}(x|\tau)\theta_{2}(y|\tau)-\theta_{1}(x|\tau)\theta_{1}(y|\tau)).
\end{gathered}
\label{eq:t2-1}
\end{equation}
\end{thm}
In the next section we first provide our proof of Theorem \ref{t2}
and then prove Theorems \ref{t1} by using Theorem \ref{t2}.

\section{Proof of Theorems \ref{t1} and \ref{t2} }

\noindent{\it Proof of Theorem \ref{t2}.} Let 
\[
F(x)=f(x)/g(x),
\]
where 
\begin{align*}
f(x) & =\theta_{1}(x+y|2\tau)\theta_{4}(x-y|2\tau)\theta_{2}(x|\tau)\theta_{2}(y|\tau)\\
 & ~-\theta_{2}(x+y|2\tau)\theta_{3}(x-y|2\tau)(\theta_{1}(x|\tau)\theta_{2}(y|\tau)+\theta_{1}(y|\tau)\theta_{2}(x|\tau))
\end{align*}
and 
\[
g(x)=\theta_{1}(x+y|2\tau)\theta_{4}(x-y|2\tau)\theta_{1}(x|\tau)\theta_{1}(y|\tau).
\]
It is easily deduced from \eqref{eq:1-4} that the entire functions
$f$ and $g$ satisfy the functional equations:
\[
h(x)=h(x+\pi)=q^{8}e^{8ix}h(x+2\pi\tau).
\]
We have $F(x)=F(x+\pi)=F(x+2\pi\tau).$ This means that $F$ is an
elliptic function with periods $\pi$ and $2\pi\tau.$ 

We temporarily assume that $0<y<\pi.$ In the fundamental period parallelogram
$\varLambda=\{a\pi+2b\pi\tau|0\leq a,b<1\},$ the function $g(x)$
has only four zeros $0,\pi\tau,\pi-y,y+\pi\tau$ and all of them are
simple. It is easily seen that 
\[
f(\pi-y)=f(y+\pi\tau)=0.
\]
 It follows from the Jacobi infinite product expressions that
\begin{align*}
2\theta_{2}(z|2\tau)\theta_{3}(z|2\tau) & =\vartheta_{2}(\tau)\theta_{2}(z|\tau),\\
2\theta_{1}(z|2\tau)\theta_{4}(z|2\tau) & =\vartheta_{2}(\tau)\theta_{1}(z|\tau).
\end{align*}
Then, by \eqref{eq:1-4} and \eqref{eq:1-3},
\[
f(0)=\vartheta_{2}(\tau)(\theta_{1}(y|2\tau)\theta_{4}(y|2\tau)\theta_{2}(y|\tau)-\theta_{2}(y|2\tau)\theta_{3}(y|2\tau)\theta_{1}(y|\tau))=0
\]
and
\[
f(\pi\tau)=q^{-2}\vartheta_{2}(\tau)(\theta_{2}(y|\tau)\theta_{1}(y|2\tau)\theta_{4}(y|2\tau)-\theta_{1}(y|\tau)\theta_{2}(y|2\tau)\theta_{3}(y|2\tau))=0.
\]
Therefore, the points $0,\pi\tau,\pi-y,y+\pi\tau$ are also zeros
of the function $f(x)$ and so the function $F(x)$ has no pole in
$\varLambda.$ Namely, $F(x)$ is an entire function of $x$ and then
$F(x)$ is a constant (independent of $x$). Set $F=C(y).$ Namely,
$f(x)=C(y)g(x).$ By analytic continuation, this identity also holds
for any complex number $y.$ That is, $F=C(y)$ holds for any complex
number $y.$ Interchanging the role of $x$ and $y$ in this identity
and noticing that $F$ is symmetric in $x$ and $y$ we can get $F=C(x).$
This means that $F$ is a constant independent of $x$ and $y,$ say,
$C.$ Namely, $F=C.$ Putting $x=y=\pi/4$ in this identity and noticing
that $\theta_{1}(\pi/4|\tau)=\theta_{2}(\pi/4|\tau)$ we can obtain
$C=1$ and then \eqref{eq:t2-1} follows. This completes the proof
of Theorem \ref{t2}. \qed

\noindent{\it Proof of Theorem \ref{t1}.} We first prove \eqref{eq:1-6}.
According to \cite[(8.7) and (8.8)]{H} we have 
\begin{equation}
\mathrm{ssn}_{q}z=\frac{\theta_{4}(z|\tau')}{\vartheta_{3}(\tau')},\:\mathrm{ccs}_{q}z=\frac{\theta_{3}(z|\tau')}{\vartheta_{3}(\tau')}.\label{eq:3-4}
\end{equation}
Replacing $\tau$ by $\tau'/2$ in \eqref{eq:t2-1}, dividing both
sides of the resulting identity by $\vartheta_{3}(\tau')\theta_{2}(x|\tau'/2)\theta_{2}(y|\tau'/2)\theta_{2}(x+y|\tau')$
and then employing \eqref{eq:0-2} and \eqref{eq:3-4} we get 
\[
\begin{gathered}\mathrm{ccs}_{q}(x-y)\mathrm{tan}_{q^{2}}x+\mathrm{ccs}_{q}(x-y)\mathrm{tan}_{q^{2}}y-\mathrm{ssn}_{q}(x-y)\mathrm{tan}_{q}(x+y)\\
=-\mathrm{ssn}_{q}(x-y)\mathrm{tan}_{q^{2}}x\mathrm{tan}_{q^{2}}y\mathrm{tan}_{q}(x+y).
\end{gathered}
\]
Then \eqref{eq:1-6} follows readily by substituting $x+y=\pi-z$
into this identity. 

The identity \eqref{eq:1-7} follows easily by dividing both sides
of \eqref{eq:1-6} by $\mathrm{tan}_{q^{2}}x\mathrm{tan}_{q^{2}}y$\\
$\cdot\mathrm{tan}_{q}z$. \qed

\section*{Acknowledgement}

This work was partially supported by the National Natural Science
Foundation of China (Grant No. 11801451).

\end{document}